\documentclass[12pt]{article}
\usepackage{latexsym}
\usepackage{amsfonts}
\usepackage{amsmath}
\usepackage{amsthm}
\usepackage{epsfig}

\newtheorem{claim}{Claim}

\newtheorem{thm}{Theorem}
\newtheorem{defn}[thm]{Definition}
\newtheorem{lemma}{Lemma}

\newtheorem{prop}[thm]{Proposition}

\newtheorem{conjecture}{Conjecture}
\newcommand{\thecons} {{16\cdot 3^d}}

\newcommand{\Z}{\mathbb Z}
\newcommand{\R}{\mathbb R}

\newcommand{\prob}{{\bf P}}

\newcommand{\be}{\begin{equation}}
\newcommand{\ee}{\end{equation}}

\def\23#1{\lfloor 2{#1}/3 \rfloor}

\def\22#1{2^{2^{#1}}}

\title{A lower bound for the chemical distance in sparse long-range percolation models}
\author{Noam Berger}

\begin{document}
\maketitle
\begin{abstract}
We consider long-range percolation in dimension $d\geq 1$, where distinct sites $x$ and $y$ are connected
with probability $p_{x,y}\in[0,1]$. Assuming that $p_{x,y}$ is translation invariant and that $p_{x,y}=\|x-y\|^{-s+o(1)}$
with $s>2d$, we show that the graph distance is at least linear with the Euclidean distance.
\end{abstract}

\section{Introduction}\label{sec:intro}
Long-range percolation (introduced by Schulman in 1983 \cite{schul}) is a percolation
model on the integer lattice $\Z^d$ in which every
two vertices can be connected by a bond.
The probability of the bond between two vertices to be open depends on
the distance between the vertices.

As it turns out, the models of most interest are those where the probability of a bond to be open decays polynomially with its length. While early works concentrated mainly on the behavior of infinite systems and in particular on critical phenomena \cite{AN,SN,B1}, more recent papers tried to understand the geometry of these graphs. A natural question in this context is what is the typical chemical distance between two points at a given Euclidean distance. This problem was first introduced in \cite{BB}, where partial answers were given. Further research was done by Coppersmith, Gamarnik and Sviridenko \cite{gamarnik} and more recently by Biskup \cite{biskup1,biskup2}. For more background and motivation, the reader is referred to \cite{BB} and \cite{biskup1}.

\subsection{The model: definitions and known results}
Let $\{p_k\}_{k\in\Z^d}$ be such that $p_k\in[0,1]$, and such that $p_k=p_{-k}$ for every $k$.
We also assume that
\begin{equation}\label{eq:srole}
0<\lim_{\|k\|\to\infty}\frac{p_k}{\|k\|^{-s}}<\infty
\end{equation}
for some $s>0$. (We also assume that the limit exists).
Let $\{\omega_{i,j}\}_{i,j\in\Z^d}$ be random variables such that $\omega_{i,j}=1$
with probability $p_{i-j}$ and $\omega_{i,j}=0$ with probability $1-p_{i-j}$, and such that $\omega_{i,j}=\omega_{j,i}$, but
otherwise the $\omega_{i,j}$-s are independent.
Consider the following graph structure on $\Z^d$: there exists an edge between $i$ and $j$ if and only if $\omega_{i,j}=1$. Let $D(x,y)$
be the (random) graph distance between $x$ and $y$, also known as the chemical distance.

The renormalization structure (see \cite{SN} and \cite{B1}) suggests the existence of five different regimes, depending on whether $s$ is smaller, greater or equal to $d$ and $2d$, and indeed in each of these regimes we see a different behavior:

For the case $s<d$, Benjamini, Kesten, Schramm and Peres proved in \cite{BKPS} (as a corollary of one of the technical lemmas) that
\[
\prob\left(D(x,y)=\left\lceil\frac{s}{d-s}\right\rceil\right)\to 1
\]
as $\|x-y\|\to\infty$.

When $s=d$, Coppersmith, Gamarnik and Sviridenko \cite{gamarnik} proved that the chemical distance scales as 
\[
\frac{\log(\|x-y\|)}{\log\log(\|x-y\|)}
\]
as $\|x-y\|\to\infty$.

When $d<s<2d$, Biskup \cite{biskup1,biskup2} proved that
\[
D(x,y)=\log(\|x-y\|)^{\Delta+o(1)}
\]
where $\Delta=\Delta(s,d)=\log(2d/s)/\log2$

When $s=2d$, very little is known. It is believed that 
\[
D(x,y)=(\|x-y\|)^{\Theta+o(1)}
\]
where $\Theta$ is some function of $d$
and the ratio
\[
\beta=\lim_{\|k\|\to\infty}\frac{p_k}{k^{-2d}}
\]
This problem is still open.

For $s>2d$, in \cite{BB} it was shown that for $d=1$, $D(x,y)$ grows linearly with $\|x-y\|$. for higher dimensions, Coppersmith, Gamarnik and Sviridenko proved in \cite{gamarnik} that 
\[
D(x,y)\geq\|x-y\|^\eta
\]
for some $0<\eta<1$ depending on $s$ and $d$.

\subsection{Main result}

The main result of the present paper is that for $s>2d$ the chemical distance scales at least linearly with the Euclidean distance, namely:

\begin{thm}\label{thm:main}
If $s>2d$ then almost surely
\begin{equation}
\liminf_{\|x\|\to\infty}\frac{D(0,x)}{\|x\|}>0.
\end{equation}
\end{thm}

The rest of the paper is devoted to the proof of Theorem \ref{thm:main} with the exception of section \ref{sec:open} that states a question still open in the $s>2d$ case.

\section{Renormalization structure}\label{sec:renorm}
In this section we describe the renormalization structure of the system. In the sections \ref{sec:incub}--\ref{sec:pfthm} we will prove Theorem \ref{thm:main}
based on this renormalization structure.

Let $\beta$ be such that
\begin{equation}\label{eq:beta}
p_k<\frac{\beta}{\|k\|^{s}}
\end{equation}
for every $k$. Let $2d<s^\prime<s$,
and let $M$ be an integer large enough so that
\begin{equation}\label{eq:defm}
100^s\beta M^{2d-s}<\frac{1}{1000\cdot2^d},
\end{equation}
\begin{equation}\label{eq:defmfac}
(Mn!)^{s-s^\prime}>n^{2s}.
\end{equation}
for every $n$,
and
\begin{equation}\label{eq:dummass}
100^s \beta M^{2d-s^\prime} (k!)^{4d-2s^\prime}<e^{-30dk}
\end{equation}
for every $k$.
Let $C_0=M$ and $C_n=n^2$ for $n>0$, and let
\[
A_n=\prod_{i=0}^{n}{C_i}=M(n!)^2.
\]

An $n$-block is the set $j+[0,A_n)^d\subseteq\Z^d$ for $j\in\Z^d$. The {\bf children} of an $n$-block $j+[0,A_n)^d$ are
the $C_n^d$ $n-1$-blocks
\[
\left\{
\left.
j+hA_{n-1}+[0,A_{n-1})^d \right|
h\in[0,C_n)^d
\right\}
\]
\begin{defn}\label{def:goodblock}
\begin{enumerate}
\item
We say that a $0$-block $Q$ is {\bf good} under the configuration $\omega$ if there is no edge of length greater that $A_0/100$ in $Q$.
\item
A $k$-block $Q$ is {\bf good} under the configuration $\omega$ if
\begin{enumerate}
\item\label{item:bigpic}
There is no edge of length greater that $A_{k-1}/100$ in $Q$,
\item\label{item:smalpic}
Among the children of $Q$, all but at most one are good, and
\item
There exists a configuration $\omega^\prime$ agreeing with $\omega$ on every pair of vertices $(x,y)$ such that at least one of $(x,y)$ is in $Q$, such that
For all $j\in\{0,+1,-1\}^d$, the block $Q+j\frac{A_{k-1}}{2}$
satisfies (\ref{item:bigpic}) and (\ref{item:smalpic}) under the configuration $\omega^\prime$.
\end{enumerate}
\end{enumerate}
\end{defn}

\begin{lemma}\label{lem:allaregood}
Let $Q_n=[0,A_n)^d$ be the $n$-block containing the origin, and let $P_n$ be the probability that $Q_n$ is not a good block. Then 
\[
\sum_{n=1}^\infty P_n < \infty
\]
\end{lemma}
\begin{proof}
We prove using induction.
By (\ref{eq:defm}),
\begin{equation}\label{eq:p0}
P_0<\frac{2^{-d}}{1000}
\end{equation}
We want to estimate the $P_k$ based on $P_{k-1}$:
\begin{eqnarray}\label{eqar:renor}
\nonumber
P_k &\leq&
\beta A_k^{2d}\cdot
\left(\frac{A_{k-1}}{100}\right)^{-s} + 2^dC_k^{2d}P_{k-1}^2\\
\nonumber
&=&
100^s \beta A_k^{2d-s}\cdot C_k^s + 2^dC_k^{2d}P_{k-1}^2\\
\nonumber
&\leq&
100^s \beta A_k^{2d-s^\prime} + 2^dC_k^{2d}P_{k-1}^2\\
\nonumber
&=&
100^s \beta M^{2d-s^\prime} (k!)^{4d-2s^\prime} + 2^dk^{4d}P_{k-1}^2\\
&\leq&
e^{-30dk} + 2^dk^{4d}P_{k-1}^2
\end{eqnarray}
where the second inequality comes from (\ref{eq:defmfac}) and the third from (\ref{eq:dummass}).
Using (\ref{eq:p0}) and (\ref{eqar:renor}) it is easy to show inductively that $P_k<2^{-d}(k+1)^{-4d}\exp(-2k)$.
\end{proof}

\section{Length of paths within blocks}\label{sec:incub}
Let $Q=[a,b)^d$ be a cube in $\Z^d$, and let $x$ and $y$ be in $C$. An {\bf $\omega$-path from $x$ to $y$ within $Q$} is
a path $x=v_1,\ldots,v_l=y$ such that $v_1,\ldots,v_l$ are all in $Q$, and the edge $(v_i,v_{i+1})$ is open under $\omega$.
The main lemma of this section is the following:
\begin{lemma}\label{lem:incub}
There exists a constant $C>0$ such that if $Q$ is a good $k$-block and $x$ and $y$ in $Q$ satisfy $\|x-y\| > A_k/2$ then
every path from $x$ to $y$ within $Q$ is of length at least $C\|x-y\|$.
\end{lemma}
\begin{proof}
We use induction to prove the following claim: There exists $C^\prime$ such that for every $k\geq \thecons$,
if $Q$ is a good $k$-block and $x$ and $y$ in $Q$ satisfy $\|x-y\| > A_k/8$ then every path from
$x$ to $y$ within $Q$ is of length at least
\begin{equation}\label{eq:indu}
C^\prime\left(\prod_{h=\thecons}^{k}{1-\frac{\thecons}{h^2}}\right)\|x-y\|.
\end{equation}
We then take
\[
C=C^\prime\left(\prod_{h=\thecons}^{\infty}{1-\frac{\thecons}{h^2}}\right)>0.
\]
To show (\ref{eq:indu}), we take $C^\prime$ to be $100/A_{\thecons-1}$. Then (\ref{eq:indu}) follows immediately for
$k=\thecons$. For the induction step, let $Q$ be a good $k$-block and $x$ and $y$ in $Q$ satisfy $\|x-y\| > A_k/8$.
Let $\gamma=(x=v_1,\ldots,v_l=y)$ be a path from $x$ to $y$ within $Q$. Then, $\|v_{i+1}-v_i\|<A_{k-1}/100$
for all $i=1,\ldots,l-1$. There exits at most one child of $Q$ that is not good, and at most one child that is not good
in any of the translation of $Q$ by elements of 
\[\left\{0,\frac{A_{k-1}}{2},-\frac{A_{k-1}}{2}\right\}^d.\] Let these not good blocks be
denoted by $B_1,B_2,\ldots,B_j\ (j\leq 3^d+1).$

Let $a_1$ be the smallest value $a$ so that $v_a\in(B_1\cup\ldots\cup B_j)$, and let $b_1$ be so that $v_{a_1}\in B_{b_1}$ (If there is more than one choice for $b_1$ we choose it arbitrarily).
Let $z_1$ be the largest value such that $v_{z_1}\in B_{b_1}$.
Inductively, let $a_{i+1}$ be the smallest value of $a$ larger than $z_i$ so that $v_a\in(B_1\cup\ldots\cup B_j)$, let $b_{i+1}$ be so that $v_{a_{i+1}}\in B_{b_{i+1}}$ and let $z_{i+1}$ be the largest value such that $v_{z_{i+1}}\in B_{b_{i+1}}$.

Let $\gamma_1=(v_1,\ldots,v_{a_1-1})$, $\gamma_2=(v_{b_1+1},\ldots,v_{a_2-1})$ and so on, up to $\gamma_n$. Let $\nu_1=(v_{a_1-1},\ldots,v_{b_1+1})$, $\nu_2=(v_{a_2-1},\ldots,v_{b_2+1})$ and so on, up to $\nu_m$. Note that both $n$ and $m$ are no larger than $3^d+1$.

For a path $\eta$, let $D(\eta)$ be the distance between its endpoints, and let $L(\eta)$ be the length of the path.

By the triangle inequality,
\[
\|x-y\|\leq
D(\gamma_1)+\ldots+D(\gamma_n)+D(\nu_1)+\ldots+D(\nu_m)
\]

Also, $D(\nu_i)<\frac{102}{100}A_{k-1}$ for every $i=1,\ldots,m$.
Let $U=\{i=1\ldots n|D(\gamma_i)>\frac{1}{2}A_{k-1}\}$. Then
\[
\sum_{i\notin U}D(\gamma_i)\leq \frac{n}{2}A_{k-1}
\]

and therefore
\begin{eqnarray}\label{eq:manygood}
\nonumber
\sum_{i\in U}D(\gamma_i)&\geq& \|x-y\| - \sum_{i=1}^mD(\nu_i) - \sum_{i\notin U}D(\gamma_i) \geq
\|x-y\| - 2(3^d)A_{k-1}\\
&\geq& \left(1-\frac{16\cdot 3^d}{k^2}\right)\|x-y\|.
\end{eqnarray}

The proof will be complete once we prove the following claim:
\begin{claim} For every $i\in U$,
\[
L(\gamma_i)\geq D(\gamma_i)\cdot C^\prime\prod_{h=\thecons}^{k-1}{1-\frac{\thecons}{h^2}}.
\]
\end{claim}
\begin{proof}
Let $v_{b_{i-1}+1}=w_1,\ldots,w_s=v_{a_{i}+1}$ be vertices in $\gamma_i$ such that for every $n$,
\begin{enumerate}
\item\label{item:big4}
$\|w_{n+1}-w_n\|>A_{k-1}/4$, and
\item\label{item:smal2}
$\|v-w_n\|<A_{k-1}/2$
for every $v\in[w_n,w_{n+1})$, where we use the notation $[w_n,w_{n+1})$ for the part of the path between $w_n$ and $w_{n+1}$.
\end{enumerate}
It is easy to see that such choice of points exists. By \ref{item:smal2}. above and the choice of $\gamma_i$, there exists a good $k-1$ block containing $[w_n,w_{n+1})$. By \ref{item:big4}. above and the induction hypothesis,
\[
L([w_n,w_{n+1}))\geq D([w_n,w_{n+1}))\cdot C^\prime\prod_{h=\thecons}^{k-1}{1-\frac{\thecons}{h^2}}
\]
and the claim follows by the triangle inequality.
\end{proof}
\end{proof}

\section{Length of inter-block paths}
\begin{prop}\label{prop:sabich}
Let $k\geq 2^d$.
There exists a constant $C>0$ such that if $Q$ is a good $k$-block
and for every $j\in\{0,+1,-1\}^d$, the $k$-block $Q+j\frac{A_k}{2}$ is good and for every $j>k$ the $j$-block
$\hat{Q}_j$ centered at $Q$ is good, then if
$x$ and $y$ in $Q$ satisfy $\|x-y\| > A_k/8$ then
every path from $x$ to $y$ is of length at least $C\|x-y\|$.
\end{prop}
\begin{proof}
There exists a block $Q_\star$ in
\[
\left\{\left.Q+j\frac{A_k}{2}\right|j\in\{0,+1,-1\}^d\right\}
\]
such that $x\in Q_\star$ and $x$ is at distance at least $A_k/4$ from the boundary of $Q_\star$. For every $j>k$,
the distance between $x$ and the boundary of $\hat{Q}_j$ is at least $A_j/4$. Let $\hat{Q}_k=Q_\star$.
Let $\gamma$ be a path between $x$ and $y$.

If $\gamma$ is included in $Q_\star$, then by Lemma \ref{lem:incub}, we are done. Otherwise, let $j_0$ be the smallest value
of $j$ such that $\gamma$ is in $\hat{Q}_j$. Let $u$ be the last point of $\gamma$ in $\hat{Q}_{j_0-1}$. Then, since $\hat{Q}_{j_0}$
is good, we get that $\mbox{dist}(u,\partial \hat{Q}_{j_0-1})<A_{j_0-1}/100$. Therefore,
$\mbox{dist}(u,x)\geq A_{j_0-1}/4$, and by Lemma \ref{lem:incub}, we are done.
\end{proof}

\section{Proof of Theorem \ref{thm:main}}\label{sec:pfthm}
Theorem \ref{thm:main} will follow easily from the following lemma:
\begin{lemma}\label{lem:dirc}
Let $v\in\R^d$ be such that $\|v\|=1$. Then
if $[nv]$ is the element of
$\Z^d$ closest to $nv$, then almost surely
\begin{equation}\label{eq:limit}
\liminf_{n\to\infty} \frac{D(0,[nv])}{\|nv\|}>0
\end{equation}
\end{lemma}

\noindent
{\bf Remark:} Note that the $D(0,[nv])$ may be infinite if $[nv]$ and $0$ are not in the same connected component.
If we force the nearest neighbor bonds to exist, then the limit exists and is finite.
If we do not force the nearest neighbor bonds to exist, then the limit on the subsequence $\{n:[nv]\mbox{ is connected to } 0 \}$
exists. We conjecture that it has to be finite, but this is not known (see Section \ref{sec:open}).

\begin{proof}
Without loss of generality, we may assume that all nearest neighbor bonds are present.
Kingman's subadditive ergodic theorem (see \cite{durrett}) guarantees that
\[
\lim_{n\to\infty} \frac{D(0,[nv])}{\|nv\|}
\]
exists almost surely. Lemma \ref{lem:allaregood} and Proposition \ref{prop:sabich} shows that it is positive.
\end{proof}
\begin{proof}[Proof of Theorem \ref{thm:main}]
Again, we assume without loss of generality that all nearest neighbor bonds are present.
By the proof of Lemma \ref{lem:dirc},
\[
C(v) = \liminf_{n\to\infty} \frac{D(0,[nv])}{\|nv\|}
\]
is bounded away from zero. Let $C=\min\{C(v)|\|v\|=1\}>0$. Let $\{v_k\}_{k=1}^{M}$ be such that for every $v$ of norm
$1$, there exists $k$ such that $\|v-v_k\|<C/2$. Then, Theorem \ref{thm:main} follows from Lemma \ref{lem:dirc} and the
triangle inequality.
\end{proof}


\section{An open problem}\label{sec:open}
In the previous sections we proved a linear lower bound for the chemical distance. We conjecture the following:
\begin{conjecture}\label{conj:uppbound}
Under the assumptions of Theorem \ref{thm:main}, almost surely,
\begin{equation}
\limsup_{\|x\|\to\infty}\frac{1_{\{0\leftrightarrow x\}}D(0,x)}{\|x\|}< \infty.
\end{equation}
\end{conjecture}
Conjecture \ref{conj:uppbound} is trivial if the nearest neighbor bonds are all present, and follows from Antal-Pisztora
under weaker assumption. However, we believe that in general, Conjecture \ref{conj:uppbound} should be hard to prove.
The following are two special cases of Conjecture \ref{conj:uppbound}, the second being a special case of the first.

\begin{conjecture}\label{conj:uppboundinf}
Under the assumptions of Theorem \ref{thm:main} with the additional assumption that
a.s. there exists an infinite cluster, almost surely,
\begin{equation}
\limsup_{\|x\|\to\infty}\frac{1_{\{0\leftrightarrow x\}}D(0,x)}{\|x\|}< \infty.
\end{equation}
\end{conjecture}

\begin{conjecture}\label{conj:uppboundinfsupcrit}
Under the assumptions of Conjecture \ref{conj:uppboundinf} with the additional assumption that
the system is super-critical, almost surely,
\begin{equation}
\limsup_{\|x\|\to\infty}\frac{1_{\{0\leftrightarrow x\}}D(0,x)}{\|x\|}< \infty.
\end{equation}
\end{conjecture}

\noindent
Noam Berger

\noindent
The California institute of Technology
\begin{verbatim}
berger@its.caltech.edu
\end{verbatim}

\end{document}